\documentclass[11pt]{article}

\usepackage{graphicx}
\usepackage{amsmath}
\usepackage{amsfonts}
\usepackage{fancyhdr} 
\usepackage{mathrsfs}
\usepackage{amssymb}
\usepackage{wrapfig}
\usepackage{bbold}
\usepackage{xcolor, soul}
\usepackage{mdframed}
\usepackage{authblk}

\newtheorem{thm}{Theorem}
\newtheorem{lem}{Lemma}

\newtheorem{prop}[lem]{Proposition}

\numberwithin{equation}{section}
\numberwithin{defn}{section}

\newcommand{\real}{\mathbb{R}}
\newcommand{\complex}{\mathbb{C}}

\newcommand{\integer}{\mathbb{Z}}
\newcommand{\ddd}{\mathbb{D}}

\newcommand{\pro}{\mathfrak{p}}

\newcommand{\ttt}{\mathbb{T}}

\newmdenv[linecolor=red,linewidth=0,leftmargin=0,rightmargin=0,leftmargin=0,rightmargin=0,skipabove=\partopsep,skipbelow=\parsep,innertopmargin=0,innerbottommargin=0,innerleftmargin=0,innerrightmargin=0]{highlight}
{\end{highlight}}
{\end{highlight}}
{\end{highlight}}

\newenvironment{pf}{\noindent {\em Proof}.\ \ }{\hspace*{\fill}\rule{.5ex}{1.4ex}\,}

\title{Inverse scattering for the one-dimensional Helmholtz equation with piecewise constant wave speed\\[5pt]
}
\begin{document}
\author[1]{Sophia Bugarija}
\author[2]{Peter C.~Gibson\footnote{Corresponding author.}}
\author[3]{Guanghui Hu}
\author[4]{Peijun Li}
\author[5]{Yue Zhao}
\affil[1]{NSERC USRA, Department~of Mathematics \& Statistics, York University, 4700 Keele St., Toronto, Ontario, Canada, M3J~1P3}
\affil[2]{Department~of Mathematics \& Statistics, York University, 4700 Keele St., Toronto, Ontario, Canada, M3J~1P3}
\affil[3]{Beijing Computational Science Research Center, Beijing 100193, People's Republic of China}
\affil[4]{Department of Mathematics, Purdue University, West Lafayette, IN, USA}
\affil[5]{School of Mathematics \& Statistics, China Central Normal University, Wuhan, Hubai 430079, People's Republic of China}

\renewcommand\Authands{ and }

\date{Version: December 9, 2019}
\makeatletter
\let\newtitle\@title
\let\newauthor\@author
\let\newdate\@date
\makeatother

\maketitle

\begin{abstract}
This paper analyzes inverse scattering for the one-dimensional Helmholtz equation in the case where the wave speed is piecewise constant.    
Scattering data recorded for an arbitrarily small interval of frequencies is shown to determine the wave speed uniquely, and a direct reconstruction algorithm is presented.  The algorithm is exact provided data is recorded for a sufficiently wide range of frequencies and the jump points of the wave speed are equally spaced with respect to travel time.  Numerical examples show that the algorithm works also in the general case of arbitrary wave speed (either with jumps or continuously varying etc.) giving progressively more accurate approximations as the range of recorded frequencies increases.   A key underlying theoretical insight is to associate scattering data to compositions of automorphisms of the unit disk, which are in turn related to orthogonal polynomials on the unit circle.  The algorithm exploits the three-term recurrence of orthogonal polynomials to reduce the required computation. 
\end{abstract}

\begin{center}MSC 35J05, 34L25, 35R05, 35R30; Keywords: one-dimensional Helmholtz equation, inverse scattering, discontinuous coefficient, reconstruction algorithm\end{center}

\newpage

\tableofcontents
\newpage

\pagestyle{fancyplain}

\section{Introduction\label{sec-introduction}}

This paper considers the one-dimensional Helmholtz equation
\begin{equation}\label{helmholtz}
u^{\prime\prime}+\kappa^2u=-\delta(x-x_0),\qquad\kappa=\omega/c,
\end{equation}
in the case where the wave speed $c(x)$ is a step function having jump points $x_j$ $(1\leq j\leq n)$, indexed according to their natural order
\[
x_0<x_1<\cdots<x_n.
\]
Thus the source point $x_0$ lies to the left of the jump points of $c$, and $c$ has the form
\begin{equation}\label{wave-speed}
c(x)=c_0\chi_{(-\infty,x_1]}(x)+\left(\sum_{j=1}^{n-1}c_j\chi_{(x_j,x_{j+1}]}(x)\right)+c_n\chi_{(x_n,\infty)}(x). 
\end{equation}
It is assumed (as usual) that radiation only travels outward from the interval $[x_0,x_n]$; this means that there is no right-travelling wave in the interval $(-\infty,x_0)$, and no left-travelling wave in $(x_n,\infty)$.  The outgoing radiation condition determines a unique $C^1$ solution to (\ref{helmholtz}).  

The goal of the paper is to analyze the inverse medium problem for (\ref{helmholtz}), as follows.   
Suppose that the source location $x_0$, measurement location $x_\ast$ situated  between $x_0$ and $x_1$, and wave speed $c_0$ at the source, are known.  To what extent can the full wave speed (\ref{wave-speed}) be determined given measured data comprised of the solution $u(x_\ast,\omega)$ to (\ref{helmholtz}) at $x_\ast$ for a range of frequencies $\omega_{\min}<\omega<\omega_{\max}$?  

\subsection{Technical preliminaries}

The standard terminology pertaining to the Helmholtz equation (left- and right-travelling waves, for example) stems from its connection to the wave equation via the Fourier transform.  This connection is detailed briefly below for the sake of clarity; afterward we fix some key notation.  

Denote the Fourier transform by
\begin{equation}\label{fourier}
\widehat{f}(\omega)=\int_\real f(t)e^{i\omega t}\,dt\quad\mbox{ where }\quad f\in L^2(\real)\cap L^1(\real),
\end{equation}
with the corresponding extension to tempered distributions as per standard theory.    
Applying the Fourier transform with respect to time to the wave equation in one spatial dimension,
\begin{equation}\label{wave}
\frac{1}{c^2}v_{tt}-v_{xx}=\delta(x-x_0)\delta(t),
\end{equation}
and writing $u(x,\omega)=\widehat{v}(x,\omega)$, leads to the Helmholtz equation (\ref{helmholtz}) with a point source at $x_0$, where $u^{\prime\prime}$ denotes the second derivative with respect to $x$ formally treating the frequency $\omega$ as a parameter rather than an independent variable.  

By the Fourier inversion theorem, any solution $u$ to (\ref{helmholtz}) yields a solution to the wave equation (\ref{wave}) of the form
\[
v(x,t)=\frac{1}{2\pi}\int_\real e^{-it\omega}u(x,\omega)\,d\omega,
\]
the integrand $e^{-it\omega}u(x,\omega)$ being referred to as a time harmonic wave of frequency $\omega$.   
Based on this association between time harmonic waves and solutions $u(x,\omega)$ to (\ref{helmholtz}), exponentials of the form
\[
e^{-i\kappa x}\quad\mbox{ and }\quad e^{i\kappa x}
\]
are referred to respectively as left-travelling and right-travelling waves, since the corresponding time harmonic waves, 
\[
e^{-it\omega}e^{-i\kappa x}=e^{-i\omega(t+x/c)}\quad\mbox{ and }\quad e^{-it\omega}e^{i\kappa x}=e^{-i\omega(t-x/c)},
\]
are left-travelling and right-travelling, respectively, with wave speed $c$.  


Fix a finite interval $\Omega=(\omega_{\min},\omega_{\max})$ with $\omega_{\min}<\omega_{\max}$, and set 
\begin{equation}\label{d}
d:\Omega\rightarrow\complex,\qquad d(\omega)=u(x_\ast,\omega)
\end{equation}
where $u$ is the unique solution to (\ref{helmholtz}) consistent with the outgoing radiation condition.  The function $d$ will be referred to as measured data.

\subsection{Overview of the paper}

The main results of the paper are to show that $d$ uniquely determines $c$ of the form (\ref{wave-speed}) no matter how small the interval of frequencies $\Omega$, and to provide an explicit, fast reconstruction algorithm.  In analyzing the inverse problem for (\ref{helmholtz}) it turns out to be more convenient to work with a surrogate reflection coefficient $R$ defined in terms of $d$ using known parameters.  In the case where $R$ is periodic the given reconstruction algorithm is exact, provided the interval $\Omega$ includes a full period.  Perhaps more importantly, for arbitrary $c$---not necessarily piecewise constant---the algorithm provides an approximate inversion method that numerical tests reveal to be remarkably accurate, the accuracy improving with increasing length of the interval $\Omega$.  A final contribution of the paper is to link the Helmholtz equation to orthogonal polynomials as well as complex analytic and almost periodic structure, providing new tools for its analysis.  

The paper is organized as follows. Connections to the existing literature are described briefly in \S\ref{sec-connection} below.  Section~\ref{sec-reformulation} introduces the surrogate reflection coefficient $R$ and derives some of its basic properties, drawing on complex analytic and almost periodic structure.  The basic uniqueness result is derived in \S\ref{sec-uniqueness}, culminating in Theorem~\ref{thm-uniqueness}, \S\ref{sec-r-delta2c}.  The fast reconstruction algorithm pertaining to the case of periodic $R$ is derived in \S\ref{sec-fast-inversion}.  Pseudocode in \S\ref{sec-fast-algorithm} summarizes the method, which was implemented in Matlab for numerical testing. Three numerical examples are presented in \S\ref{sec-numerical}. The paper ends with a brief discussion in \S\ref{sec-discussion}.

\subsection{Connection to known results\label{sec-connection}}

Recent work on inverse scattering for the one-dimensional Helmholtz equation has centred on iterative methods based either on the Volterra scattering series with both reflection and transmission data \cite{Ya:2014}, or on optimization using a convex cost functional \cite{Kl1:2016,Kl2:2017,Kl3:2018}, requiring that the wave speed be $C^2$.  The case of piecewise constant wave speed considered in the present paper has until now been open.  Iterative reconstruction methods that have been proposed for the one-dimensional Helmholtz tend to be computationally expensive, raising the question of whether more computationally efficient direct methods are possible.  

Direct reconstruction methods have been developed for the one-dimensional wave equation with piecewise constant coefficients in \cite{Gi:JCP2018}.  The present paper relates to the latter work, demonstrating that the link to orthogonal polynomials can be exploited also in the case of the Helmholtz equation.  The translation between the wave and Helmholtz equations is not immediate, however.   The Fourier transform of the Helmholtz equation with a source term, as considered in the present paper, leads to a wave equation with a $\delta$ forcing term, the corresponding wavefield being a regular function comprised of a superposition of travelling square waves.  By contrast, the wavefield in \cite{Gi:JCP2018} is the convolution of a  purely singular superposition of Dirac functions with a smooth source wave form.  

Orthogonal polynomials on the unit circle may be represented as products of $2\times 2$ matrices \cite{SiOPUC1:2005}, which is what relates them to compositions of automorphisms of the unit disk.  A backward recurrence formula for the Fourier transform of the impulse response for the wave equation with piecewise constant wave speed, which involves disk automorphisms, was known long before the connection to OPUC came to light. Extended to the torus, this backward recurrence formula has a Fourier series whose coefficients were described completely in \cite{Gi:JFAA2017}; a simple consequence of that description is used in the present paper.

\section{Reformulation of the data\label{sec-reformulation}}

The present section introduces a surrogate reflection coefficient $R(\omega)$ defined in terms of $u(x_\ast,\omega)$, and derives its basic properties.  The importance of $R$ is that on one hand it is determined by the measured data $d$, as will be shown below, while on the other hand it readily exhibits certain complex analytic and almost periodic structure that is not obviously a feature of the measured data itself.  These deeper mathematical structures are key to solving the inverse medium problem for (\ref{helmholtz}).

\subsection{Definition and basic structure of $R$}

Fix notation as follows. Let $\pro$ denote the projective map
\begin{equation}\label{projection}
\pro:\complex^2\setminus\{0\}\rightarrow\complex\cup\{\infty\},\qquad\pro\binom{a}{b}=\frac{a}{b}.
\end{equation}
Given a matrix of the form $M=\begin{pmatrix}\zeta&\rho\zeta\\ \overline{\rho}&1\end{pmatrix}$, where $\zeta,\rho\in\complex$, and letting $v$ denote a complex variable, set
\begin{equation}\label{phi-M}
\varphi_M(v)=\pro\left(M\binom{v}{1}\right)=\zeta\frac{v+\rho}{1+\overline{\rho}v}.
\end{equation}
Thus if $\zeta\in S^1$ lies on the unit circle and $\rho\in\ddd$ is in the open unit disk then $\varphi_M$ is a disk automorphism; in general it is a linear fractional transformation.  A useful elementary fact from complex analysis (used freely in what follows) is that composition of linear fractional transformations corresponds to multiplication of matrices: if 
$M^\prime=\begin{pmatrix}\zeta^\prime&\rho^\prime\zeta^\prime\\ \rule{0pt}{12pt}\overline{\rho^\prime}&1\end{pmatrix}$ is a second matrix, then
\begin{equation}\label{composition-rule}
\varphi_{MM^\prime}=\varphi_M\circ\varphi_{M^\prime}. 
\end{equation}
Referring to $c$ of the form (\ref{wave-speed}), set
\begin{equation}\label{delta-rj-zj}
\Delta_j=\frac{x_j-x_{j-1}}{c_{j-1}},\quad
r_j=\frac{1-\frac{c_{j-1}}{c_j}}{1+\frac{c_{j-1}}{c_j}},\quad z_j=e^{2i\omega\Delta_j}\qquad(1\leq j\leq n)
\end{equation}
and write 
\begin{equation}\label{Mj}
M_j=\begin{pmatrix}z_j&r_jz_j\\ \bar{r}_j&1\end{pmatrix}\qquad(1\leq j\leq n).
\end{equation}
Note that $-1<r_j<1$ and $\varphi_{M_j}$ is a disk automorphism for each $1\leq j\leq n$.  Because $r_j\in\real$, the complex conjugation in the $(2,1)$-entry of $M_j$ plays no role and can be ignored.  However, it will be seen later that the same matrices $M_j$ arise in the context of orthogonal polynomials, in which case the $r_j\in\ddd$ may be complex valued; the (redundant) complex notation is included in (\ref{Mj}) to make clearer the connection.

Now comes the key definition of this section: given $c_0,x_0,x_\ast$ and $u(x_\ast,\omega)$, set 
\begin{equation}\label{R-defintion}
R(\omega)=-\frac{2ic_0}{\omega}e^{i\frac{\omega}{c_0}(x_\ast-x_0)}u(x_\ast,\omega)-e^{2i\frac{\omega}{c_0}(x_\ast-x_0)}.
\end{equation}
\begin{prop}\label{prop-matrix-product}
\[
R(\omega)=\pro\left(M_1\cdots M_n\binom{0}{1}\right)
=\varphi_{M_1}\circ\cdots\circ\varphi_{M_n}(0)
\]
\end{prop}
\begin{pf}
The second equality follows from (\ref{phi-M}) and (\ref{composition-rule}).  The first equality follows from direct analysis of the Helmholtz equation (\ref{helmholtz}) with coefficient of the form (\ref{wave-speed}) as follows. 

For $x<x_1$ the general solution to (\ref{helmholtz}) compatible with the outgoing radiation condition works out to be
\begin{equation}\label{left-solution}
u(x,\omega)=\left\{
\begin{array}{cc}
\left(A_0-\frac{\kappa}{2i}\right)e^{-i\kappa(x-x_0)}&\mbox{ if }x<x_0\\
\rule{0pt}{18pt}A_0e^{-i\kappa(x-x_0)}-\frac{\kappa}{2i}e^{i\kappa(x-x_0)}&\mbox{ if }x_0<x<x_1
\end{array}\right.,
\end{equation}
where $A_0\in\complex$ is arbitrary.  In the intervals $x_j<x<x_{j+1}$ $(1\leq j\leq n-1)$ and $x>x_n$, equation (\ref{helmholtz}) reduces to the homogeneous equation 
\begin{equation}\label{homogeneous}
u^{\prime\prime}+\kappa^2u=0
\end{equation}
which has general solution of the form 
\begin{equation}\label{homogeneous-solution}
u(x,\omega)=
\left\{\begin{array}{cc}
A_je^{-i\kappa(x-x_j)}+B_je^{i\kappa(x-x_j)}&\mbox{ if }x_j<x<x_{j+1}\quad(1\leq j\leq n-1)\\
\rule{0pt}{18pt}B_ne^{i\kappa(x-x_n)}&\mbox{ if }x>x_n
\end{array}\right..
\end{equation}
Note that the outgoing radiation condition prohibits any left-moving component of the form $A_ne^{-i\kappa(x-x_n)}$ in the interval $x>x_n$.  

Fixing $B_0=-\frac{\kappa}{2i}$ (as prescribed by (\ref{left-solution})) and $A_n=0$, the requirement $u\in C^1$ (with respect to $x$) imposes relations among the various pairs $\binom{A_j}{B_j}$ as follows.  Continuity of $u$ at $x_j$ requires
\begin{equation}\label{continuity-0}
A_{j-1}e^{-i\omega\Delta_j}+B_{j-1}e^{i\omega\Delta_j}=A_j+B_j\qquad(1\leq j\leq n).
\end{equation}
Continuity of $u^\prime$ at $x_j$ forces
\begin{equation}\label{continuity-1}
\frac{-A_{j-1}}{c_{j-1}}e^{-i\omega\Delta_j}+\frac{B_{j-1}}{c_{j-1}}e^{i\omega\Delta_j}=\frac{-A_j}{c_j}+\frac{B_j}{c_j}\qquad(1\leq j\leq n). 
\end{equation}
The continuity equations (\ref{continuity-0}) and (\ref{continuity-1}) combine to yield the fundamental relation
\begin{equation}\label{fundamental-relation}
\binom{A_{j-1}}{B_{j-1}}=\frac{e^{-i\omega\Delta_j}(c_j+c_{j-1})}{2c_j}M_j\binom{A_j}{B_j}\qquad(1\leq j\leq n),
\end{equation}
where $M_j$ is defined in accordance with (\ref{delta-rj-zj}) and (\ref{Mj}). 

It follows from (\ref{Mj}) and the prescribed values of $B_0=-\frac{\kappa}{2i}$ and $A_n=0$ that 
\begin{equation}\label{revised-matrix}
\binom{A_0}{-\frac{\kappa}{2i}}=\gamma M_1\cdots M_n\binom{0}{B_n}\quad\mbox{ where }\quad\gamma=\prod_{j=1}^n\frac{e^{-i\omega\Delta_j}(c_j+c_{j-1})}{2c_j}.
\end{equation}
The above equation guarantees $B_n\neq0$ if $\omega\neq0$, since $-\frac{\kappa}{2i}\neq0$.  For $\omega\neq0$, applying the projective map (\ref{projection}) to (\ref{revised-matrix}) yields the relation
\begin{equation}\label{after-projection}
\frac{A_0}{B_0}=-\frac{2i}{\kappa}A_0=\pro\left(M_1\cdots M_n\binom{0}{1}\right).  
\end{equation}
According to (\ref{left-solution}), since $x_0<x_\ast<x_1$, $u(x_\ast,\omega)$ has the form
\[
u(x_\ast,\omega)=A_0e^{-i\frac{\omega}{c_0}(x_\ast-x_0)}-\frac{\omega}{2ic_0}e^{i\frac{\omega}{c_0}(x_\ast-x_0)}.
\] 
Therefore $A_0/B_0$ may be expressed in terms of $u(x_\ast,\omega)$ and the quantities $c_0, x_0$ and $x_\ast$ as
\begin{equation}\label{coefficient-data}
\frac{A_0}{B_0}=-\frac{2ic_0}{\omega}e^{i\frac{\omega}{c_0}(x_\ast-x_0)}u(x_\ast,\omega)-e^{2i\frac{\omega}{c_0}(x_\ast-x_0)}=R(\omega).
\end{equation}
The desired result then follows from (\ref{after-projection}). 
\end{pf}

\subsection{Further properties of $R$}

Instead of fixing $z_j=e^{2i\omega\Delta_j}$ in the definition of $M_j$ as in (\ref{delta-rj-zj}), (\ref{Mj}), one can alternatively let $(z_1,\ldots,z_n)\in\overline{\ddd}^n$ vary over the $n$-dimensional (closed) polydisk, thereby extending $R$ to a function
\begin{equation}\label{f-defn}
f:\overline{\ddd}^n\rightarrow\complex,\qquad f(z_1,\ldots,z_n)=\varphi_{M_1}\circ\cdots\circ\varphi_{M_n}(0),
\end{equation}
so that $R$ is in turn the restriction of $f$ to the torus line
\begin{equation}\label{torus-line}
\omega\mapsto \ell(\omega)=(e^{2i\omega\Delta_1},\ldots,e^{2i\omega\Delta_n}).
\end{equation}
Consideration of $f$ allows one to infer various properties of $R$ relatively easily.  
\begin{prop}\label{prop-f-properties}
\begin{enumerate}
\item For every $(z_1,\ldots,z_n)\in\overline{\ddd}^n$
\[
|f(z_1,\ldots,z_n)|\leq \tanh\left(\sum_{j=1}^n\tanh^{-1}|r_j|\right)<1.
\]
\item 
$f:\overline{\ddd}^n\rightarrow \ddd$ is holomorphic.
\item $R=f\circ\ell:\real\rightarrow\ddd$ is holomorphic.
\end{enumerate}
\end{prop}
\begin{pf}
For $r_0\in\ddd$ and $0\leq \varepsilon\leq1$, it is easily verified that
\[
\max_{v\in\varepsilon\overline{\ddd}}\left|\frac{v+r_0}{1+\overline{r_0}v}\right|=\frac{\varepsilon+|r_0|}{1+\varepsilon|r_0|}=\tanh\left(\tanh^{-1}\varepsilon+\tanh^{-1}|r_0|\right).
\]
Applying this iteratively to the composition $f(z_1,\ldots,z_n)=\varphi_{M_1}\circ\cdots\circ\varphi_{M_n}(0)$ yields 
\[
\max_{z_j\in\overline{\ddd}}\left|f(z_1,\ldots,z_n)\right|\leq \tanh\left(\tanh^{-1}|r_1|+\cdots+\tanh^{-1}|r_n|\right),
\]
proving 1.  

To prove 2., consider a single variable $z_j$, with the other variables being fixed.  Denote by $w_j$ the constant
\[
w_j=\frac{\varphi_{M_{j+1}\cdots M_n}(0)+r_j}{1+\bar{r}_j\varphi_{M_{j+1}\cdots M_n}(0)}
\]
so that $|w_j|<1$ by 1.  Thus for $z_j\in\overline{\ddd}$, $z_jw_j\in\ddd$, and so holomorphicity of the linear fractional transformation $\varphi_{M_1\cdots M_{j-1}}$ on $\ddd$ implies that 
\begin{equation}\label{f-single-variable}
f(z_1,\ldots,z_n)=\varphi_{M_1\cdots M_{j-1}}(z_jw_j)
\end{equation}
is holomorphic with respect to $z_j$ in $\frac{1}{|w_j|}\ddd\supset\overline{\ddd}$.  It follows by Hartogs' theorem that $f$ is jointly holomorphic in $z_1,\ldots,z_n$ on the closed polydisk $\overline{\ddd}^n$.  

The mapping (\ref{torus-line}) is the restriction to $\real$ of an entire function; combined with 2., this proves 3. 
\end{pf}

Proposition~\ref{prop-f-properties} implies that for any $r=(r_1,\ldots,r_n)\in\ddd^n$, the restriction of $f$ to the $n$-torus $\ttt^n$ is a $C^\infty$ function whose Fourier series converges pointwise (see \cite[Ch.~3]{Gr:2008}).  The Fourier coefficients of $f$ are functions of $r$ and have a highly non-trivial structure first described in \cite{Gi:JFAA2017}, the main result of which is as follows. 

For each $(p,q)\in\integer^2$ write $\psi^{(p,q)}:\complex\rightarrow\complex$ to denote the polynomial defined as follows. If $\min\{p,q\}\geq 1$ set
\begin{equation}\label{varphi}
\psi^{(p,q)}(\zeta)=\textstyle\frac{(-1)^p}{q(p+q-1)!}\,\displaystyle(1-\zeta\bar{\zeta})\frac{\partial^{\,p+q}}{\partial\bar{\zeta}^p\partial\zeta^q}(1-\zeta\bar{\zeta})^{p+q-1}.
\end{equation}
If $\min\{p,q\}<0$ or $p=0<q$ set $\psi^{(p,q)}=0$; and if $p\geq 0$ set $\psi^{(p,0)}(\zeta)=\zeta^p$. Note in particular that $\psi^{(0,0)}(\zeta)=1$.   
The following theorem adopts the convention that indices beyond $n$ should be interpreted as zero, so that for example, if $k\in\integer^n$, then $k_{n+1}=0$.   
\begin{thm}[See {\cite[Thm.~1]{Gi:JFAA2017}}]\label{thm-fourier-expansion} For $z=(z_1,\ldots,z_n)\in\ttt^n$ and fixed $r=(r_1,\ldots,r_n)\in\ddd^n$, 
the Fourier series of the $n$-fold composition of disk automorphisms $f(z)=\varphi_{M_1}\circ\cdots\circ\varphi_{M_n}$ is
\begin{equation}\label{fourier-expansion}
f(z)=\sum_{k\in\{1\}\times\integer_+^{n-1}}\biggl(\,\prod_{j=1}^n\psi^{(k_j,k_{j+1})}(r_j)\biggr)z^k.
\end{equation}
\end{thm}

The full force of the above theorem is not required for present purposes---just the corollary comprising part 3 of the following result. 
\begin{prop}\label{prop-ap}
Suppose without loss of generality that $r_1\neq 0$.  Then $R$ may be uniquely expressed in the form 
\[
R(\omega)=\sum_{j=1}^\infty a_je^{it_j\omega},
\]
where $t_1<t_2<\cdots$ and each $a_j\neq 0$.  Furthermore:
\begin{enumerate}
\item $\sum_{j=1}^\infty|a_j|^2<1$;
\item $R$ is almost periodic in the sense of Besicovitch; 
\item $t_1=2\Delta_1$ and $a_1=r_1$.  
\end{enumerate}
\end{prop}
\begin{pf}
Replacing $z$ in (\ref{fourier-expansion}) with $\ell(\omega)=(e^{2i\omega\Delta_1},\ldots,e^{2i\omega\Delta_n})$ and collecting terms with common frequency yields the given form for $R$.  Holomorphicity of $f$ on $\overline{\ddd}^n$ implies absolute summability of the coefficients, from which it follows that 
\[
\sum_{j=1}^\infty |a_j|^2<\infty.
\]
It follows in turn that $R$ is almost periodic in the sense of Besicovitch (see \cite[p.~505]{Be:1926}).  Therefore, by the Plancherel theorem for Besicovitch almost periodic functions,
\[
\begin{split}
\sum_{j=1}^\infty |a_j|^2&=\lim_{L\rightarrow\infty}\frac{1}{2L}\int_{-L}^L|R(\omega)|^2\,d\omega\\
&\leq \tanh^2\left(\sum_{j=1}^n\tanh^{-1}|r_j|\right)<1,
\end{split}
\]
where the inequality follows from part 1~of Proposition~\ref{prop-f-properties}.  Finally, the lowest frequency in (\ref{fourier-expansion}) that occurs with non-zero coefficient after substituting $z=\ell(\omega)$  is $t_1=2\Delta_1$, corresponding to the multi-index $k=(1,0,\ldots,0)$.  Since $\psi^{(0,0)}(r_j)=1$, the corresponding amplitude simplifies to 
$\psi^{(1,0)}(r_1)=r_1$.  
\end{pf}

The foregoing results allow one to draw on complex analysis and the theory of almost periodic functions to prove that the measured data $d$ uniquely determines the coefficient (\ref{wave-speed}).

\section{Uniqueness and reconstruction\label{sec-uniqueness}}

\subsection{$d\mapsto R$\label{sec-d2R}}

\begin{prop}\label{prop-d2R}
The measured data $d$ determines $R$.
\end{prop}
\begin{pf}
Recall that the measured data $d$ is the restriction of $u(x_\ast,\omega)$ to $\Omega=(\omega_{\min},\omega_{\max})$. Therefore the restriction of $R$ to $\Omega$ may be expressed in terms of $d$ as 
\[
R|_\Omega(\omega)=-\frac{2ic_0}{\omega}e^{i\frac{\omega}{c_0}(x_\ast-x_0)}d(\omega)-e^{2i\frac{\omega}{c_0}(x_\ast-x_0)}.
\]
Part 3 of Proposition~\ref{prop-f-properties} implies $R$ may be recovered from $R|_\Omega$ by analytic continuation.  
\end{pf}

It should be noted, however, that analytic continuation is numerically impractical, a fact first demonstrated in \cite{Sa:1989}.  A practical alternative is discussed below in $\S\ref{sec-fast-inversion}$.  

\subsection{$R\mapsto (r,\Delta)$\label{sec-R2r-delta}}

The sequences of reflectivities $r=(r_1,\ldots,r_n)$ and layer widths $\Delta=(\Delta_1,\ldots,\Delta_n)$ are determined by $R$ as follows.  Suppose that the index $n$ in the representation (\ref{wave-speed}) is minimal, so that each of the $r_j$ defined in (\ref{delta-rj-zj}) is non-zero.  Using Proposition~\ref{prop-ap}, write 
\[
R(\omega)=r_1e^{2i\Delta_1\omega}+\sum_{j=2}^\infty a_je^{it_j\omega},
\]
with $t_j>2\Delta_j$ for all $j\geq 2$.  Then, according to standard theory (see \cite{Be:1926}), 
\begin{equation}\label{ap-delta1}
2\Delta_1=\inf\left\{\lambda\,\left|\,\lim_{L\rightarrow\infty}\frac{1}{2L}\int_{-L}^LR(\omega)e^{-i\lambda\omega}\,d\omega\neq 0\right.\right\}
\end{equation}
and 
\begin{equation}\label{ap-r1}
r_1=\lim_{L\rightarrow\infty}\frac{1}{2L}\int_{-L}^LR(\omega)e^{-2i\Delta_1\omega}\,d\omega.
\end{equation}

Now, each of the disk automorphisms $\varphi_{M_j}$ is invertible, with inverse
\begin{equation}\label{inverse}
\varphi_{M_j}^{-1}(v)=\varphi_{M_j^{-1}}(v)=\bar{z_j}\frac{v-r_jz_j}{1-\overline{r_jz_j}v}.  
\end{equation}
Thus the pair $(r_1,\Delta_1)$ determines
\begin{equation}\label{layer-stripping}
\varphi_{M_1}^{-1}\bigl(R(\omega)\bigr)=\varphi_{M_2}\circ\cdots\circ\varphi_{M_n}(0),
\end{equation}
effectively decreasing the number of layers by one. 
Iterating this process determines successive pairs $(r_2,\Delta_2),\ldots,(r_n,\Delta_n)$, terminating in the zero function.  

The formulas (\ref{ap-delta1}) and (\ref{ap-r1}) suffice theoretically to determine $\Delta_1$ and $r_1$ (and hence the sequences $\Delta$ and $r$ by iteration), but neither has an efficient numerical implementation.  Convergence of the mean integral may be slow, making it difficult in practice to distinguish zero from non-zero limits.  On the other hand, given the pair of sequences $(r,\Delta)$, it is a straightforward matter to determine the function $c$, both theoretically and numerically. 

\subsection{$(r,\Delta)\mapsto c$\label{sec-r-delta2c}}

Given $x_0$ and $c_0$, the pair of sequences $(r,\Delta)$ determines the the wave speed $c$ as follows.  Note first that by (\ref{delta-rj-zj}),
\begin{equation}\label{rj-reformulation}
r_j=\frac{c_j-c_{j-1}}{c_j+c_{j-1}}=\tanh\log\sqrt\frac{c_j}{c_{j-1}}\qquad(1\leq j\leq n).
\end{equation}
Therefore the values $c_1,\ldots,c_n$ may be expressed in terms of $r$ and $c_0$ by the formula
\begin{equation}\label{r2c}
c_j=c_0\exp\left({2\sum_{i=1}^j\tanh^{-1}r_i}\right)\qquad(1\leq j\leq n). 
\end{equation}
The next step is to reconstruct the jump points $x_j$ of $c$ in terms of $\Delta$ and the values $c_j$. It follows directly from the formula for $\Delta_j$ in (\ref{delta-rj-zj}) that 
\begin{equation}\label{c2x}
x_j=x_0+\sum_{i=1}^jc_{i-1}\Delta_i\qquad(1\leq j\leq n),
\end{equation}
allowing one to reconstruct the wave speed
\[
c(x)=c_0\chi_{(-\infty,x_1]}(x)+\left(\sum_{j=1}^{n-1}c_j\chi_{(x_j,x_{j+1}]}(x)\right)+c_n\chi_{(x_n,\infty)}(x)
\]
as desired.  To summarize, 
\begin{thm}\label{thm-uniqueness}
Fix $\omega_{\min}<\omega_{\max}$ and let $u$ denote the unique solution to (\ref{helmholtz}) consistent with the outgoing radiation condition, and such that the wave speed $c$ has the form (\ref{wave-speed}). 
Measured data of the form 
$
d(\omega)=u(x_\ast,\omega)
$
for all $\omega_{\min}<\omega<\omega_{\max}$
uniquely determines $c$.  
\end{thm}

\section{Fast reconstruction for equal layers\label{sec-fast-inversion}}

As noted earlier, the theoretical reconstruction of $c$ from $d$ outlined in \S\ref{sec-uniqueness}
does not readily translate into a practical algorithm.  There is an alternate approach, however, that does allow for efficient numerical implementation if the layer thicknesses $\Delta_j$ are all the same, provided data is collected for a sufficiently wide bandwidth of frequencies. It turns out that this latter class of equal-layer-thickness media (in which $c_j=c_{j+1}$ is allowed) can serve to approximate essentially arbitrary wave speed profiles $c$---not just step functions---leading to an efficient inversion algorithm with wide applicability.  

\subsection{Equal layer thicknesses, periodic $R$\label{sec-equal-layer}}

The case of equal layer thicknesses relates to the theory of orthogonal polynomials on the unit circle, an association that first came to light in the context of the wave equation \cite{Gi:NMPDE2018,Gi:JCP2018}.  The present case of the Helmholtz equation exhibits similar underlying structure---although the associated wave equation (\ref{wave}) is different from that in \cite{Gi:NMPDE2018,Gi:JCP2018}.  

A preliminary step is to consider the function $f$ defined in (\ref{f-defn}) in the special case of equal layers. 
Suppose there is a fixed value $\Delta_0$ such that 
\begin{equation}\label{constant-hypothesis}
\Delta_j=\Delta_0\qquad(1\leq j\leq n).
\end{equation}
Setting $z=e^{2i\Delta_0\omega}$, the notation (\ref{delta-rj-zj}) implies that 
\begin{equation}\label{constant-z}
z_j=z\qquad(1\leq j\leq n). 
\end{equation}
Letting $z$ vary over the closed unit disk $\overline{\ddd}$, the function $f$ defined according to (\ref{f-defn}) restricts to its diagonal value
\begin{equation}\label{defn-g}
g(z)=f(z,\ldots,z)=\varphi_{M_1}\circ\cdots\circ\varphi_{M_n}(0).
\end{equation}
Note that, while the right-hand side of (\ref{defn-g}) is formally the same as before, in the present case each of the matrices $M_j$ involves the same variable $z_j=z$.  Thus 
$
g:\overline{\ddd}\rightarrow\ddd
$
is holomorphic by part 2 of Proposition~\ref{prop-f-properties}; moreover $g(0)=0$ by (\ref{Mj}). It follows that the function
\begin{equation}\label{F-defn}
F(z)=\frac{1+g(z)}{1-g(z)}:\overline{\ddd}\rightarrow\complex_+
\end{equation}
is holomorphic, has values with positive real part and takes the value $F(0)=1$.  Note that $F$ may be defined equivalently as follows. Set
\begin{equation}\label{M0}
M_0=\begin{pmatrix}1&1\\ -1&1\end{pmatrix}.
\end{equation}
Then 
\begin{equation}\label{F-defn-2}
\begin{split}
F(z)&=\varphi_{M_0}\circ\varphi_{M_1}\circ\cdots\circ\varphi_{M_n}(0)\\
&=\pro\left(M_0\cdots M_n\binom{0}{1}\right)
\end{split}
\end{equation}

The function $R(\omega)$ is obtained from $g$ by restricting to $z=e^{2i\Delta_0\omega}$,
\begin{equation}\label{periodic-R}
R(\omega)=g\negthinspace\left(e^{2i\Delta_0\omega}\right). 
\end{equation}
The Fourier coefficients of $R$ are thus the Taylor coefficients of $g$; i.e., 
\begin{equation}\label{fourier-taylor}
R(\omega)=\sum_{j=1}^\infty \alpha_je^{2ij\Delta_0\omega}\quad\mbox{ and }\quad g(z)=\sum_{j=1}^\infty \alpha_j z^j.
\end{equation}
And $R$ is periodic with period $p=\pi/\Delta_0$.  

The Herglotz representation theorem asserts that for any function of the form (\ref{F-defn}) with $F(0)=1$, there exists a probability measure $d\mu$ on the unit circle $S^1$ such that 
\begin{equation}\label{herglotz}
F(z)=\int_{S^1}\frac{\zeta+z}{\zeta-z}\,d\mu(\zeta)\qquad(z\in\ddd).
\end{equation}
Moreover, it follows directly from the formula (\ref{herglotz}) that the Taylor coefficients of $F$ are the conjugate moments of $d\mu$.  More precisely,
\begin{equation}\label{F-moments}
F(z)=1+2\sum_{j=1}^\infty m_jz^j\quad\mbox{ where }\quad m_j=\int_{S^1}\bar{\zeta}^j\,d\mu(\zeta). 
\end{equation}

\subsection{Fourier coefficients $\mapsto$ moments}

Observe by (\ref{fourier-taylor}) and (\ref{F-defn}) that the first $j$ Fourier coefficients $(\alpha_1,\ldots,\alpha_j)$ of $R$ determine the moments $(m_1,\ldots,m_j)$ of $d\mu$, for any $j\geq 1$, as follows.  Since $d\mu$ is a probability measure, the zeroeth moment $m_0=1$ is fixed. 
Referring to (\ref{fourier-taylor}) and (\ref{F-moments}), fix $j\geq 1$ and let $A$ denote the $(j+1)\times(j+1)$ lower triangular Toeplitz matrix whose first column is $(0,\alpha_1,\alpha_2,\ldots,\alpha_j)^t$.  Straightforward linear algebraic manipulations yield the equation
\begin{equation}\label{alpha2moments}
(I-A)\begin{pmatrix}1\\ m_1\\ \vdots\\ m_j\end{pmatrix}=\begin{pmatrix}1\\ 0\\ \vdots\\ 0\end{pmatrix}\end{equation}
which may be efficiently solved for the moments of $d\mu$ in terms of the Fourier coefficients of $R$ by back substitution.  In particular, 
\begin{equation}\label{m-recurrence}
m_j=\alpha_jm_0+\alpha_{j-1}m_1+\cdots+\alpha_1m_{j-1}.
\end{equation}
Since $m_0=1$ is determined, the recurrence (\ref{m-recurrence}) allows successive determination of the moments of $d\mu$ given the Fourier coefficients of $R$, with $m_j$ being determined by the $j$ Fourier coefficients $\alpha_1,\ldots,\alpha_j$, as claimed.

\subsection{Orthogonal polynomials on the unit circle\label{sec-OPUC}}

The following notation is chosen to be consistent with 
the encyclopedic survey \cite{SiOPUC1:2005,SiOPUC2:2005} of orthogonal polynomials on the unit circle (OPUC)---Szeg\H{o}'s book \cite{Sz:1975} being the classical reference.  Some standard facts will also be summarized, proofs of which can be found in \cite{SiOPUC1:2005}. 

Fix a probability measure $\mu$ on the unit circle (supported on an infinite set), and denote by
\[
\Phi_0(z)=1, \Phi_1(z), \Phi_2(z),\ldots
\]
the sequence of monic polynomials obtained by Gram-Schmidt orthogonalization of the monomial sequence $1, z, z^2,\ldots$ with respect to the $L^2(d\mu)$ inner product
\begin{equation}\label{inner-product}
\langle h,k\rangle_{d\mu}=\int_{S^1}h(z)\overline{k(z)}\,d\mu(z).
\end{equation}
For each $j\geq 0$, let $\Phi_j^\ast$ denote the dual polynomial 
\begin{equation}\label{dual-polynomial}
\Phi_j^\ast(z)=z^j\overline{\Phi_j(1/\bar{z})}.
\end{equation}
For OPUC the classical three-term recurrence takes the following form. 
\begin{prop}[{See \cite[Thm.~11.4.2]{Sz:1975}}]\label{prop-recurrence}
There exists a uniquely determined sequence of scalars $r_j\in\ddd$ $(1\leq j<\infty)$, such that 
\begin{equation}\label{3-term-recurrence}
\Phi_{j+1}(z)=z\Phi_j(z)-\overline{r_{j+1}}\Phi_j^\ast(z)\qquad(0\leq j<\infty).
\end{equation}
\end{prop}
A useful fact following from (\ref{dual-polynomial}), (\ref{inner-product}) and the fact that the $\Phi_j$ are monic is that 
\begin{equation}\label{useful-fact}
\langle \Phi_j^\ast,1\rangle_{d\mu}=\langle\Phi_j,\Phi_j\rangle_{d\mu}\neq 0,
\end{equation}
provided $\mu$ is supported on an infinite set. 
Observe that the sequence $(r_j)_{j=1}^n$ determines $(\Phi_j)_{j=0}^n$ by (\ref{3-term-recurrence}). Denote by $(\Psi_j)_{j=1}^n$ the sequence of monic orthogonal polynomials determined the sequence of recurrence coefficients $(-r_j)_{j=1}^n$.  (One can of course take $n=\infty$ in the foregoing statements.)

With the given notation in hand and exploiting the recurrence (\ref{3-term-recurrence}), one may verify by induction that 
\begin{equation}\label{OPUC-matrices}
M_0M_1\cdots M_n=\begin{pmatrix}\Psi_n(z)&\Psi_n^\ast(z)\\ -\Phi_n(z)&\Phi_n^\ast(z)\end{pmatrix}
\end{equation}
where the matrices $M_j$ have precisely the form (\ref{M0}) and (\ref{Mj}), but with $z$ in place of $z_j$.  
Thus the function $g(z)$ has an explicit representation in terms of orthogonal polynomials,
\begin{equation}\label{g-representation}
g(z)=\varphi_{M_1\cdots M_n}(0)=\varphi_{M_0}^{-1}\circ\varphi_{M_0M_1\cdots M_n}(0)=\frac{\Psi_n^\ast(z)-\Phi_n^\ast(z)}{\Psi_n^\ast(z)+\Phi_n^\ast(z)}.
\end{equation}
The crucial fact is that the polynomials $\Phi_j$ $(0\leq j\leq n)$ are orthogonal with respect to the measure $\mu$ determined by $F$ via the Herglotz representation theorem.
\begin{thm}[From {\cite[Thm.~1]{Gi:JCP2018}}]\label{thm-orthogonality}
Let $n\geq 1$, fix a sequence $(r_j)_{j=1}^n$ of points in $\ddd$, define $F$ by (\ref{F-defn-2}), and let $\mu$ be the probability measure associated to $F$ by the Herglotz representation theorem (\ref{herglotz}). For $0\leq j\leq n$, set
\[
M_0M_1\cdots M_j=\begin{pmatrix}\Psi_j(z)&\Psi_j^\ast(z)\\ -\Phi_j(z)&\Phi_j^\ast(z)\end{pmatrix}.
\]
Then the sequence $(\Phi_j)_{j=0}^n$ is orthogonal in $L^2(d\mu)$.  
\end{thm}
 It follows directly from the matrix formulation (\ref{OPUC-matrices}) that the recurrence coefficients for the polynomials $(\Phi_j)_{j=0}^n$ are precisely the reflectivities $r_1,\ldots,r_n$.  The moments of $d\mu$ allow one to calculate scalar products of polynomials, and hence to compute the reflectivities $r_j$, starting with the Fourier coefficients $\alpha_j$ of $R$ using the three-term recurrence (\ref{3-term-recurrence}).  
 If the interval $\Omega=(\omega_{\min},\omega_{\max})$ contains a full period $p=\pi/\Delta_0$ the Fourier coefficients of $R$ may be computed by the formula
\begin{equation}\label{fourier-coefficients}
\alpha_j=\frac{\Delta_0}{\pi}\int_{\omega_{\min}}^{\omega_{\min}+p}R(\omega)e^{-2ij\Delta_0\,\omega}\,d\omega.
\end{equation}
This gives rise to an algorithm to recover $c$ from $d$ as follows. 

\newpage 

\subsection{A fast inversion algorithm\label{sec-fast-algorithm}}

\begin{table}[th]
\rule{0.75\oddsidemargin}{0in}\begin{tabular}{|l|l|}
\hline
\multicolumn{2}{|c|}{\rule[-5pt]{0pt}{17pt}\textbf{Computation of $d\mapsto c$}}\\
\hline
\emph{Parameters}&\rule{0pt}{12pt}$c_0,x_0,x_\ast,\Omega=(\omega_{\min},\omega_{\max})$ such that\\[0pt]
& $c_0> 0$, $x_0<x_\ast$ and $\omega_{\min}<\omega_{\max}$\\[5pt]
\emph{Input}&\rule{0pt}{12pt}Measured data $d(\omega)\quad (\omega\in\Omega)$\\[5pt]
\emph{Preliminary step}&Set \\[0pt]
& $R(\omega)=-\frac{2ic_0}{\omega}e^{i\frac{\omega}{c_0}(x_\ast-x_0)}d(\omega)-e^{2i\frac{\omega}{c_0}(x_\ast-x_0)}\quad(\omega\in\Omega)$\\[5pt]
&Define \\[0pt]
& $p=\inf\left\{0<y<|\Omega|\,\left| \, \omega,\omega+y\in\Omega\Rightarrow R(\omega)=R(\omega+y)\right.\right\}$\\[5pt]
& Choose an integer $n\geq 1$\\[10pt]
\emph{Step 1: $d\mapsto \alpha$} & Set $\alpha_j=\displaystyle\frac{1}{p}\int_{\omega_{\min}}^{\omega_{\min}+p}R(\omega)e^{-2ij\frac{\pi}{p}\omega}\,d\omega\qquad(1\leq j\leq n)$\\[15pt] 
\emph{Step 2: $\alpha\mapsto m$} & Set $m_0=1$\\
&For $ j=1:n$, set\\[5pt]
&$\rule{18pt}{0pt}m_j=\alpha_j m_0+\alpha_{j-1}m_1+\cdots+\alpha_1 m_{j-1}$\\[10pt]
\emph{Step 3: $m\mapsto r$} & Set $\nu_0^0=1$\\
& For $j=0:n-1$, set\\[5pt]
& $\rule{18pt}{0pt}r_{j+1}={\sum_{i=0}^j\overline{\nu_i^jm_{i+1}}}\left/{\sum_{i=0}^j \nu_i^jm_{j-i}}\right.$\\[5pt]
& $\rule{18pt}{0pt}
\bigl(\nu_0^{j+1},\ldots,\nu_{j+1}^{j+1}\bigr)=\bigl(0,\nu_0^j,\ldots,\nu_j^j\bigr)-\overline{r_{j+1}}\left(\overline{\nu_j^j},\overline{\nu_{j-1}^j},\ldots,\overline{\nu_0^j},0\right)$\\[10pt]
\emph{Step 4: $r\mapsto c$} & Set $c_j=c_0\exp\left(2\sum_{i=1}^j\tanh^{-1}r_i\right)\qquad(1\leq j\leq n)$\\[5pt]
& and $x_j=x_0+\frac{\pi}{p}\sum_{i=1}^jc_{i-1}\qquad(1\leq j\leq n)$\\[10pt]
\emph{Output}&$c(x)=c_0\chi_{(-\infty,x_1]}(x)+\left(\sum_{j=1}^{n-1}c_j\chi_{(x_j,x_{j+1}]}(x)\right)+c_n\chi_{(x_n,\infty)}(x)$\\[10pt]
\hline
\end{tabular}\hspace*{\fill}
\end{table}

Note that the integer $n$ chosen in the preliminary step must be large enough for $c$ to be properly reconstructed. More precisely, $n$ determines a travel-time distance $s=n\pi/p$, and the algorithm reconstructs the restriction of $c$ to the interval $(x_0,x_0+S)$, where the time it takes a pulse to propagate from $x_0$ to $x_0+S$ is $s$.  Thus the larger $p$ is, the larger $n$ must be chosen in order to recover the restriction of $c$ to a given interval.  

\newpage

\section{Numerical examples\label{sec-numerical}}

In numerical tests involving step functions $c$ having equal layer thicknesses the above algorithm gives essentially perfect reconstruction, even with 10\% noise added, as illustrated below in \S\ref{sec-equal}.  Perhaps more interesting is the fact that the case of equal layer thickness can be used as an approximate reconstruction scheme for data coming either from a step function having unequal layer thicknesses or a continuously varying $c$.  For the latter cases, illustrated in \S\ref{sec-unequal} and \S\ref{sec-continuously} below, one sets $p=\omega_{\max}-\omega_{\min}$ artificially instead of determining it from the measured data.  The quality of the reconstructions is such that they are difficult to distinguish visually from the original wave speed; to compensate, the relative $L_2$ error of the reconstructions is recorded in the captions. 

Each of the following three examples includes two reconstructions, one with noise-free data and a second where the data has 10\% i.i.d.~Gaussian noise added.  In working with noisy data, it turns out that the reconstruction is much more stable if one evaluates the coefficients $\alpha_j$ of $R$ in Step~1 using an interval $\Omega$ that is shifted away from 0.  The shift used for noisy data in each of the three cases is recorded in the captions. In all of the plots $x_0=0$, $x_\ast=(x_1+x_0)/2$, and $N$ indicates the number of sample points used to compute the integral formulation of $\alpha_j$ by a Riemann sum.  Both clean and noisy data are plotted over the first, unshifted interval $\Omega$ to allow comparison between the two. 

\subsection{Equal layer thicknesses\label{sec-equal}}

\begin{figure}[h]
\fbox{
\includegraphics[width=1.45in]{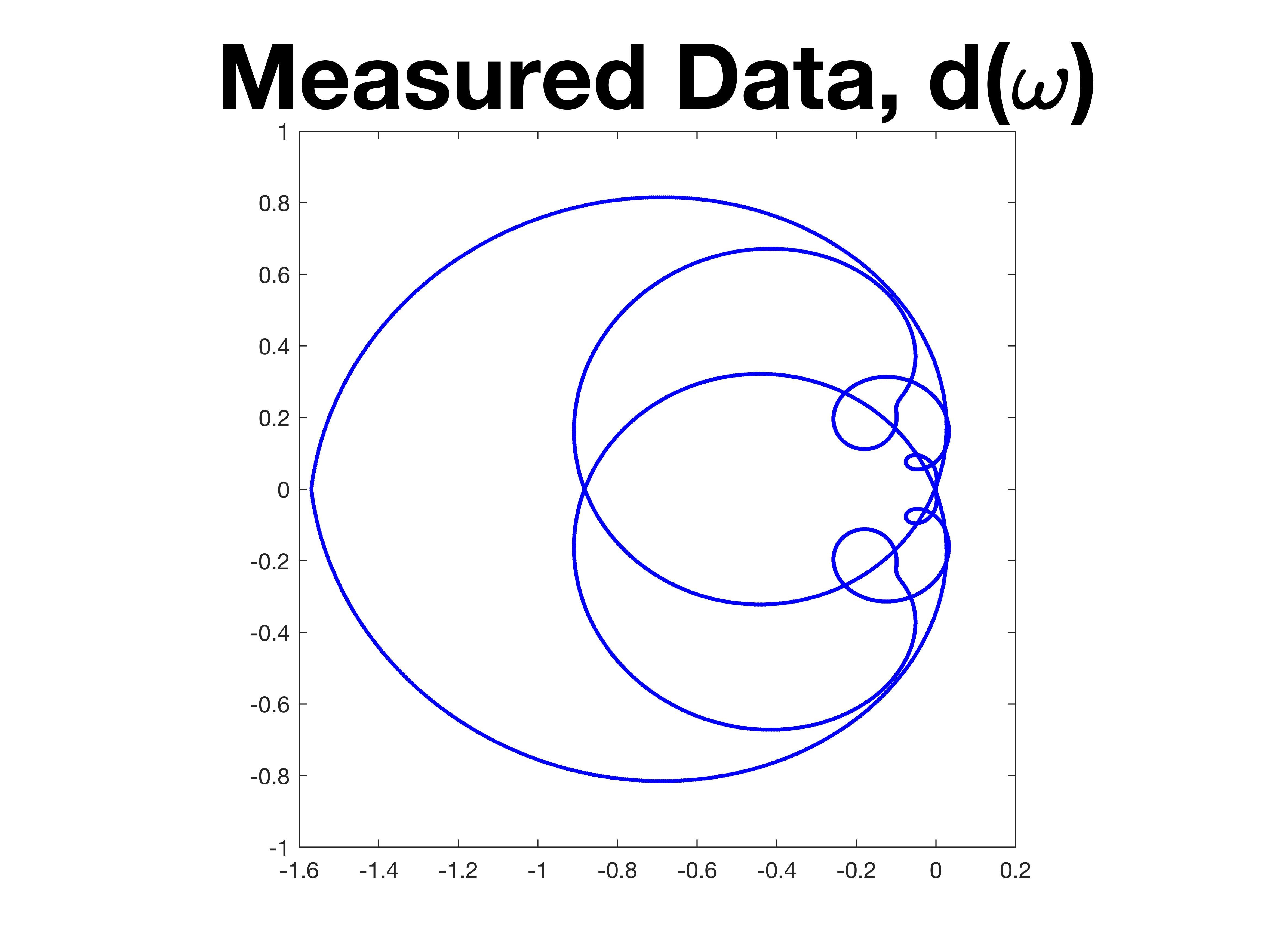}
\includegraphics[width=1.45in]{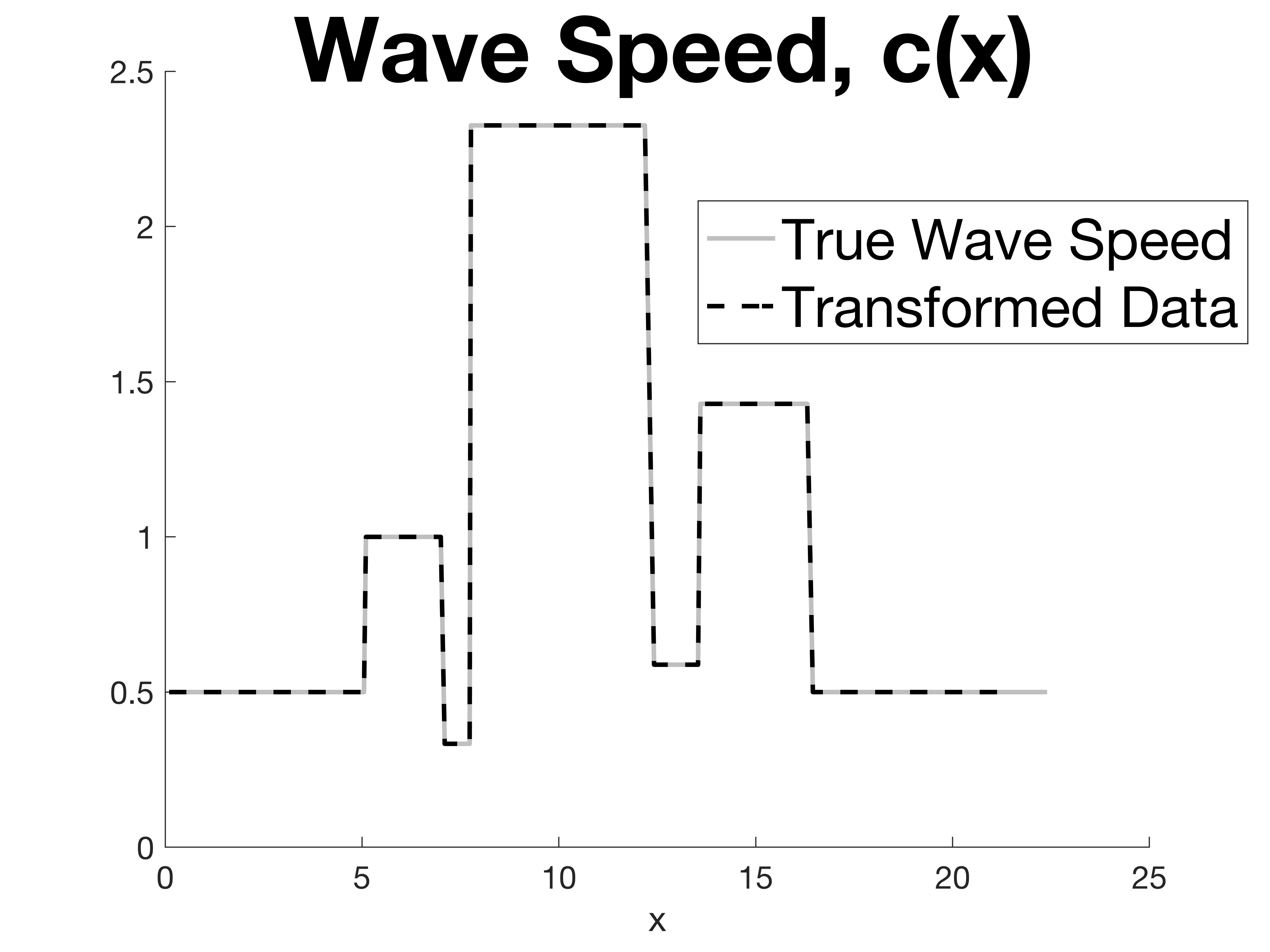}
\includegraphics[width=1.45in]{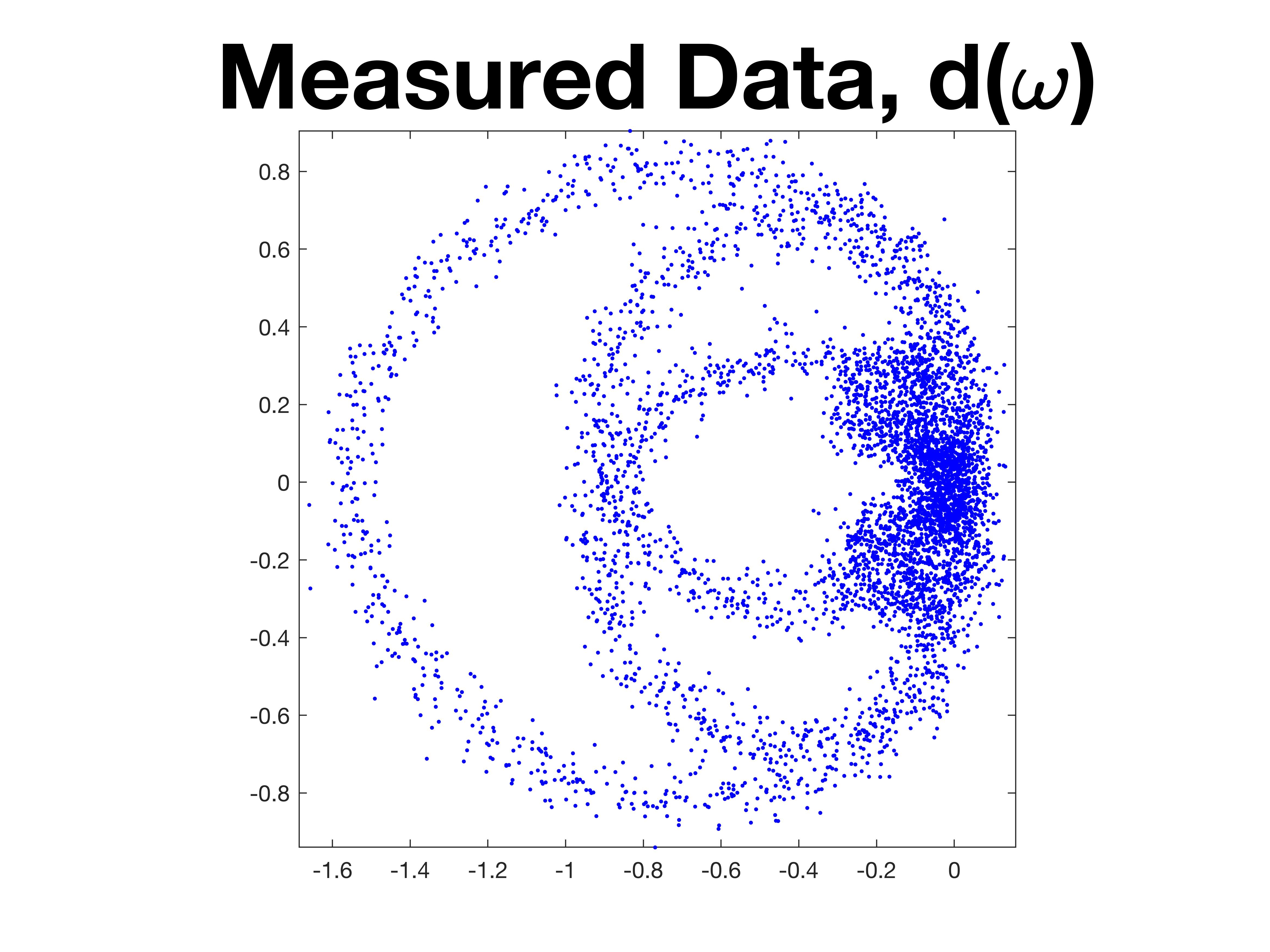}
\includegraphics[width=1.45in]{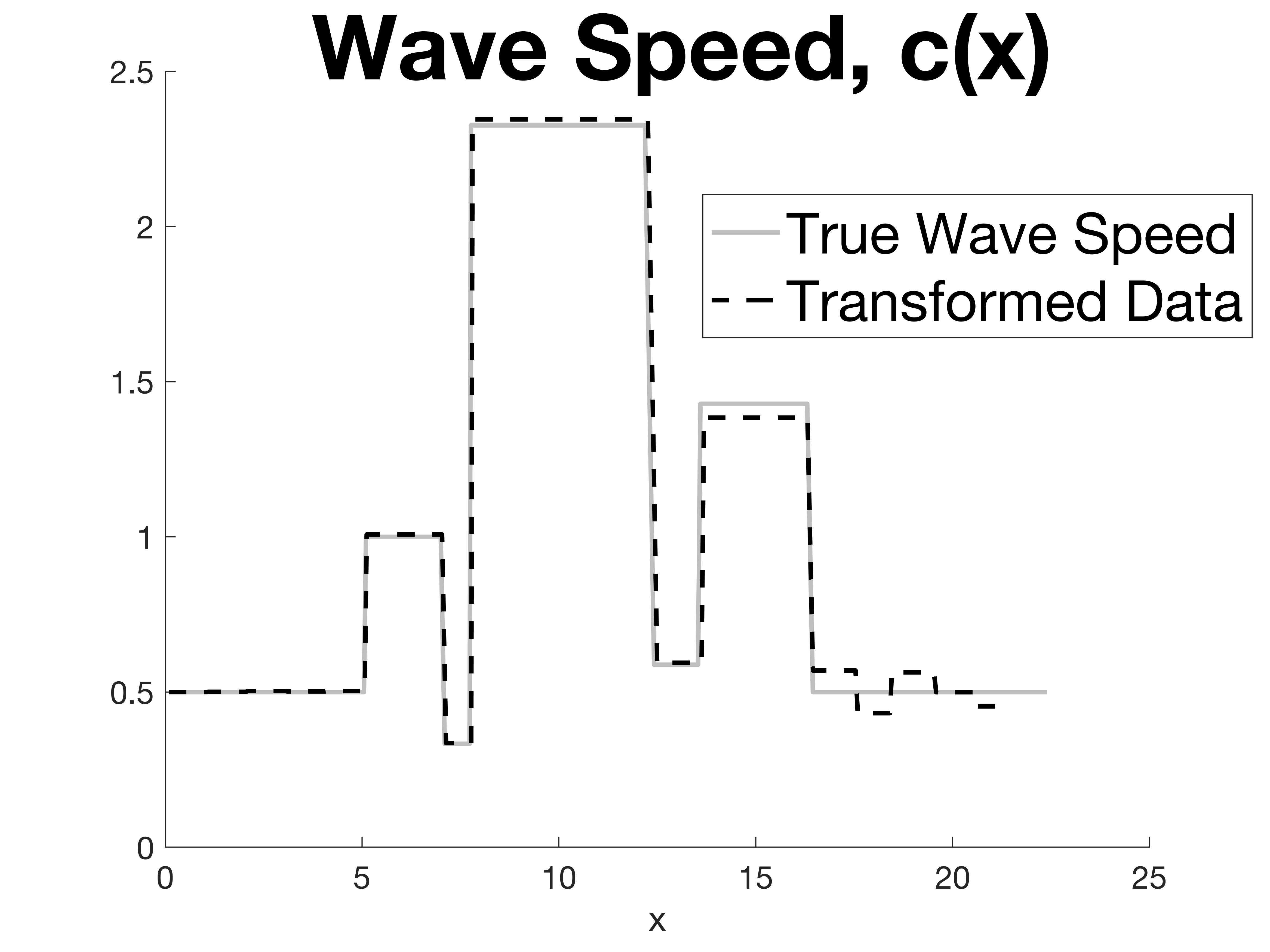}
}
\caption{A step function with jumps equally spaced with respect to travel time. $c_0=\frac{1}{2},p=\pi/2, N=5000, n=15$. At left:  $\Omega=(-.7852, .7852)$, relative error$=2.3413e-14$, computing time $1.72$sec. At right, with 10\% i.i.d.~noise:  $\Omega=(399.2148, 400.7852)$, relative error$=0.0473$, computing time $1.72$sec.}
\end{figure}

\newpage

\subsection{Unequal layer thicknesses\label{sec-unequal}}

\begin{figure}[h]
\fbox{
\includegraphics[width=1.45in]{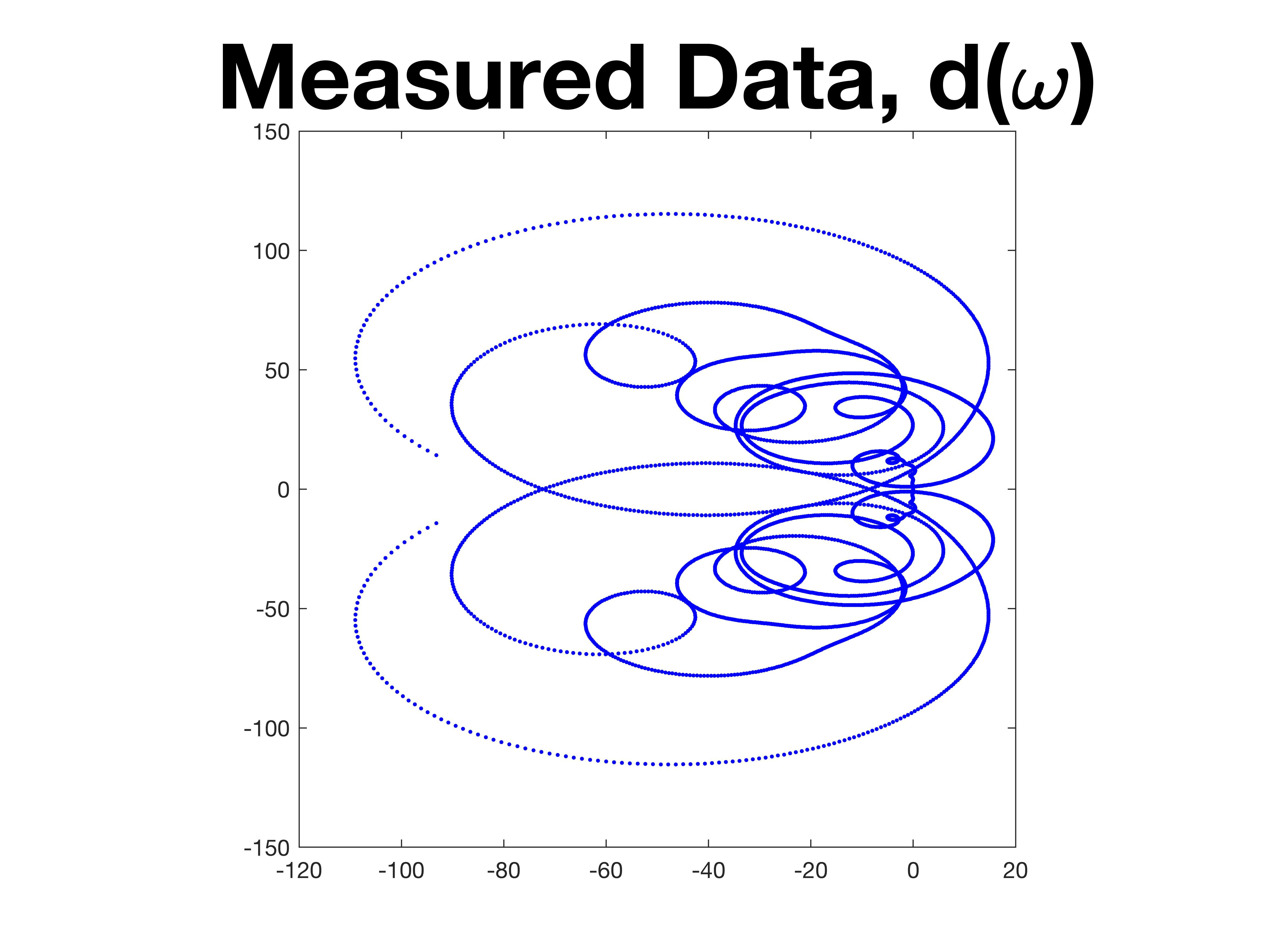}
\includegraphics[width=1.45in]{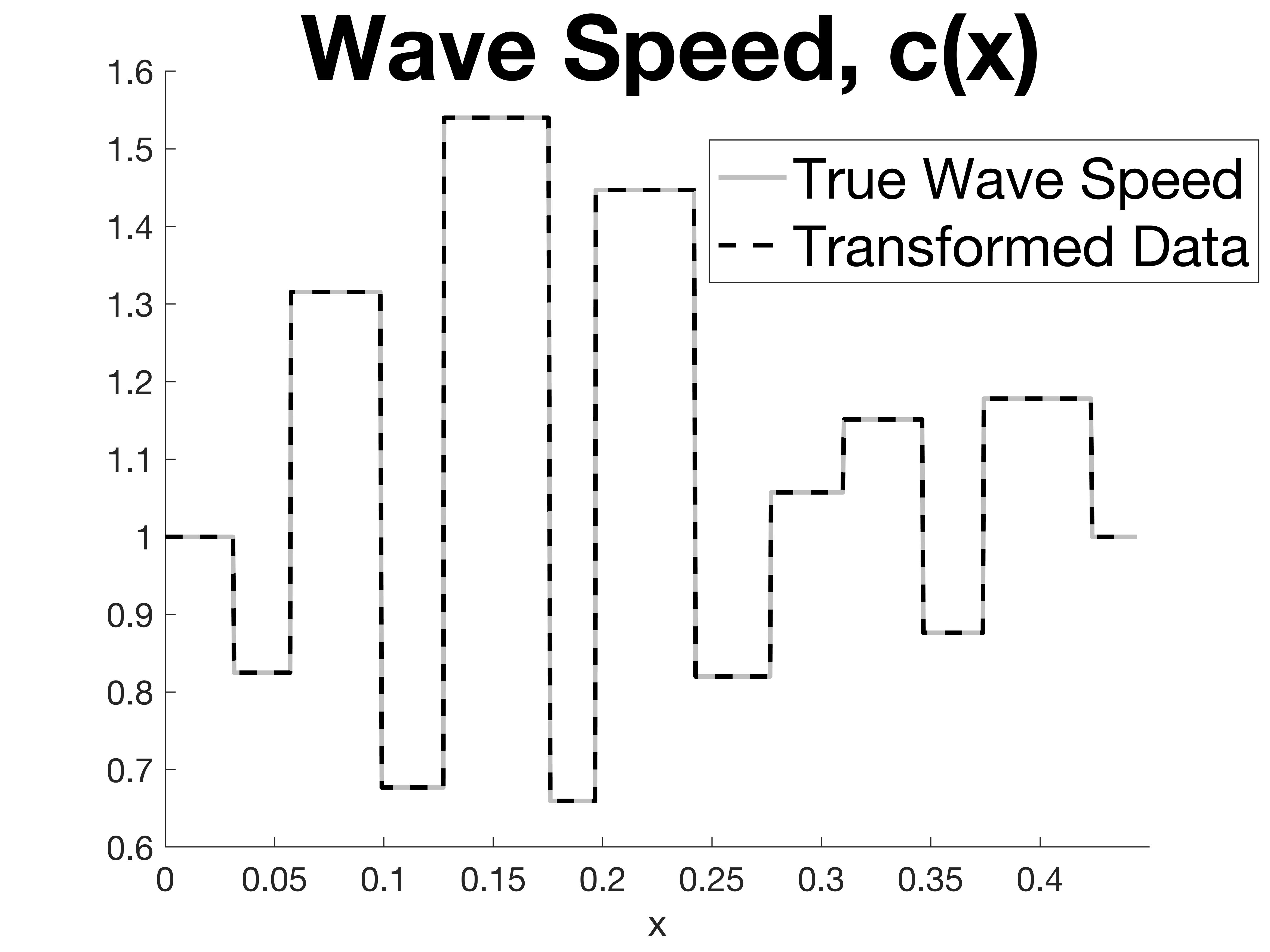}
\includegraphics[width=1.45in]{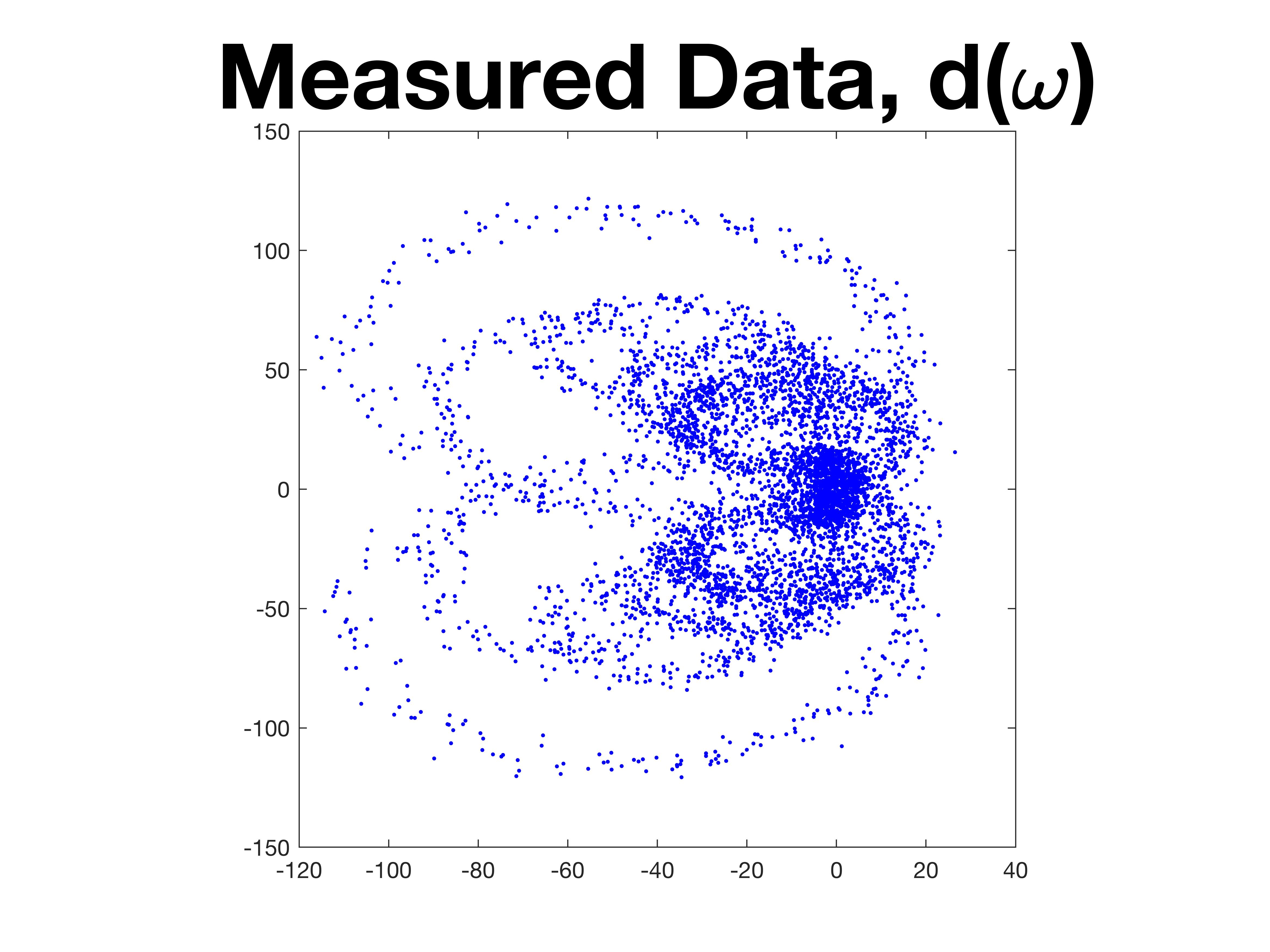}
\includegraphics[width=1.45in]{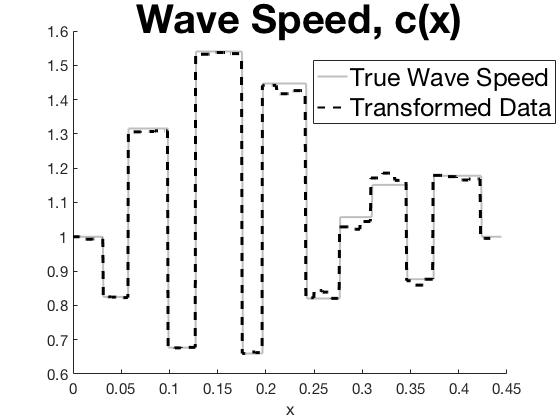}
}
\caption{A step function with jumps irregularly spaced with respect to travel time. $c_0=1, p=\pi/.0105, N=5000, n=40$. At left:  $\Omega=(-149.3919, 149.3919)$, relative error$=3.224e-15$, computing time $2.80$sec. At right, with 10\% i.i.d.~noise:  $\Omega=(250.6081, 549.3919)$, relative error$=0.0125$, computing time $1.95$sec.}
\end{figure}

\subsection{Continuously varying wave speed\label{sec-continuously}}

\begin{figure}[h]
\fbox{
\includegraphics[width=1.45in]{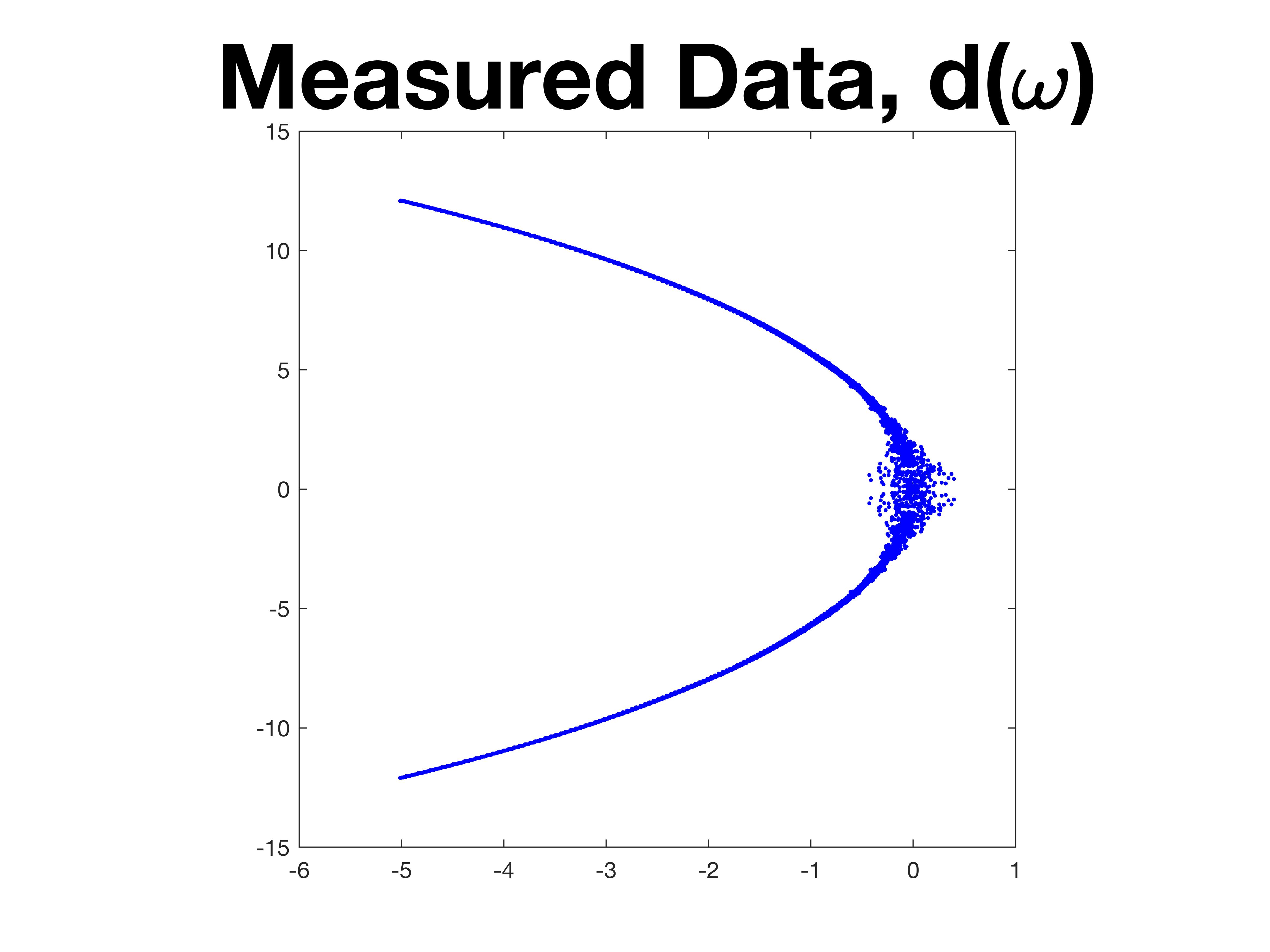}
\includegraphics[width=1.45in]{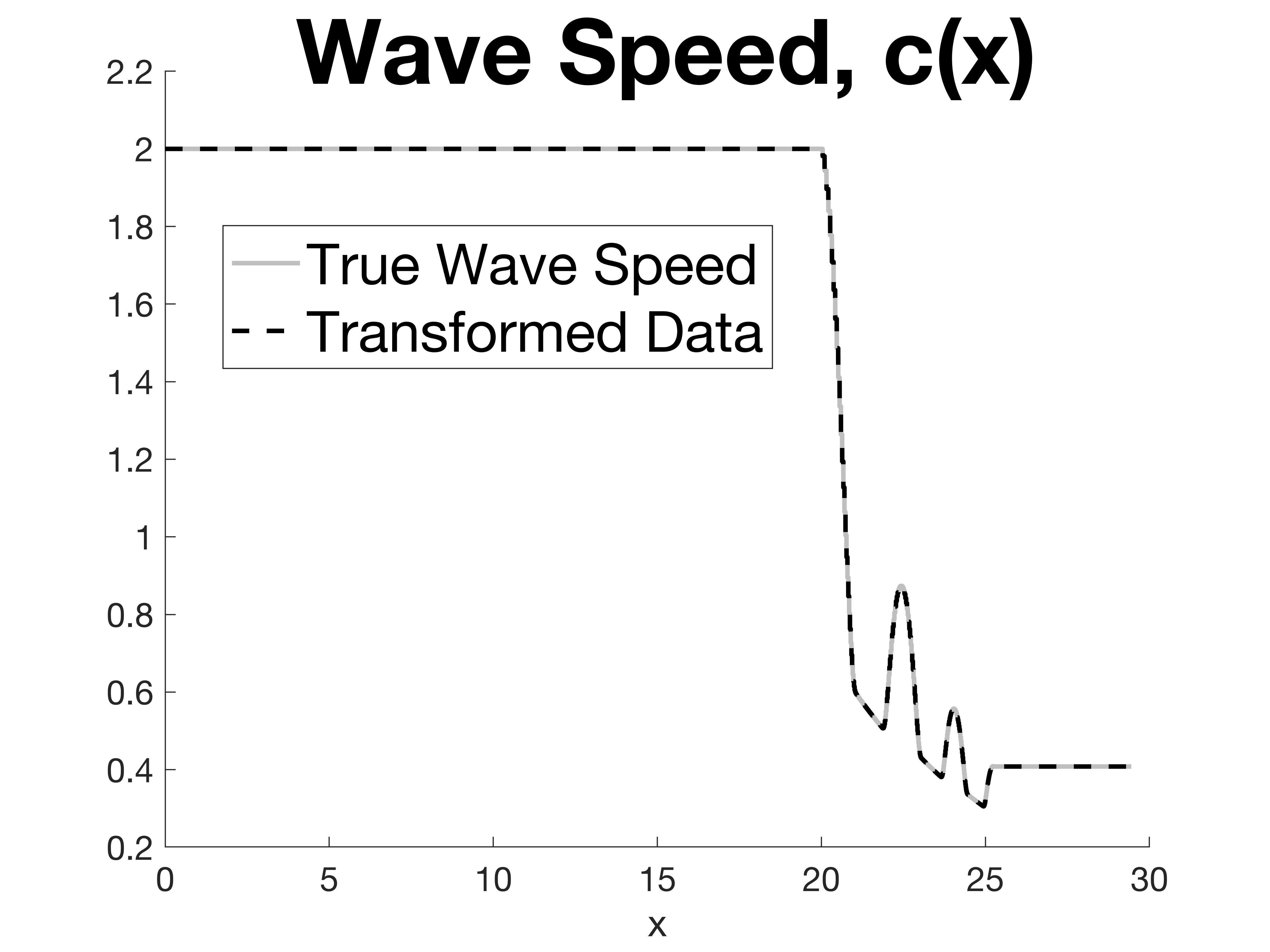}
\includegraphics[width=1.45in]{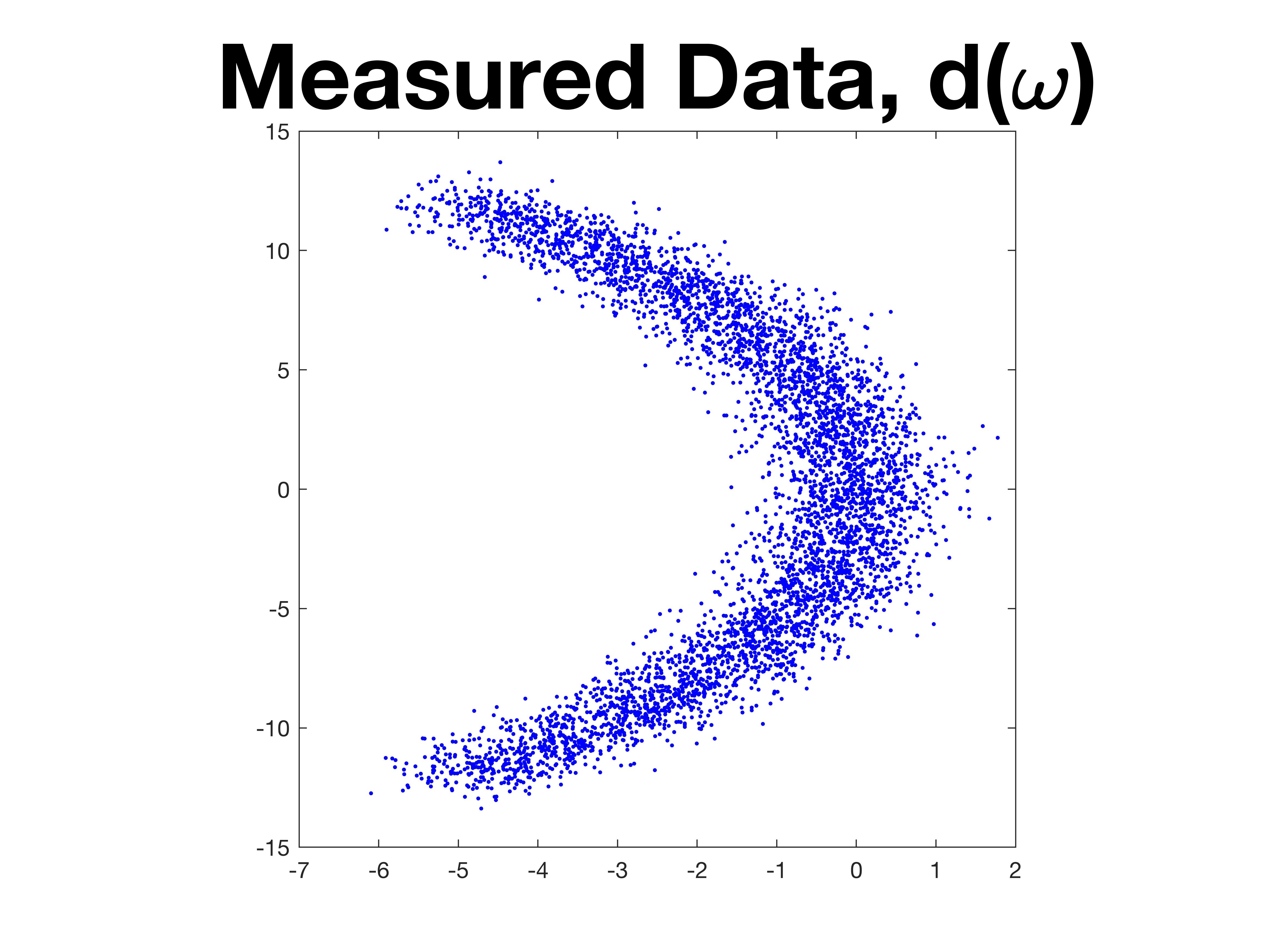}
\includegraphics[width=1.45in]{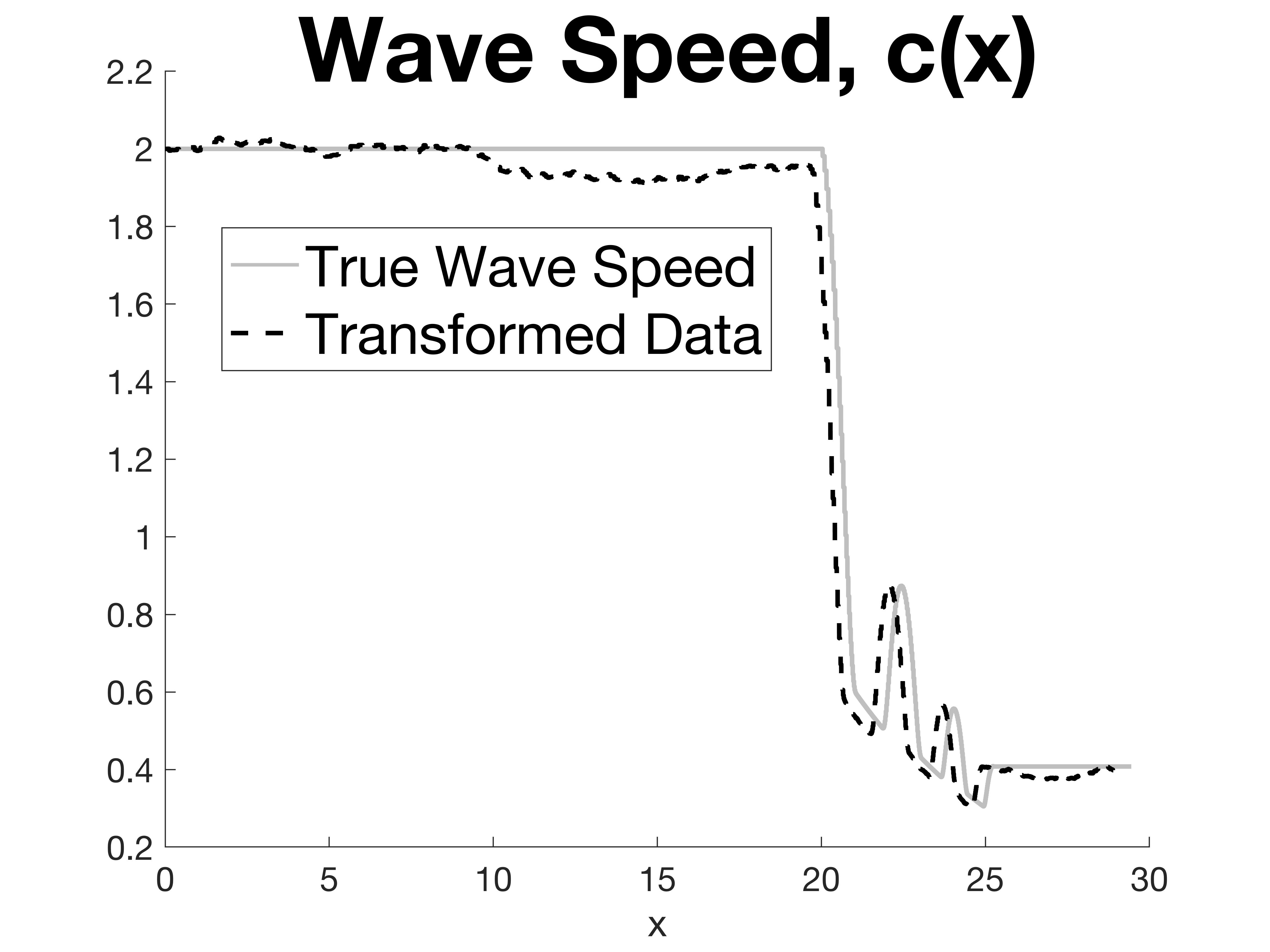}
}
\caption{A continuous function. $c_0=2,p=\pi/.03,N=5000, n=1000$. At left:   $\Omega=(-52.3555, 52.3555)$, relative error$=7.3035e-15$, computing time $36.1$sec. At right, with 10\% i.i.d.~noise:   $\Omega=(347.6506, 452.3494)$, relative error$=0.0254$, computing time $34.7$sec.}
\end{figure}

\section{Discussion\label{sec-discussion}}

Theorem~\ref{thm-uniqueness} and the explicit algorithm described in \S\ref{sec-fast-algorithm} are the two principal contributions of the present paper, both grounded in the analytic framework established in \S\ref{sec-reformulation}.  The speed of the inversion algorithm depends crucially on the connection to OPUC as formulated in \S\ref{sec-equal-layer} and \S\ref{sec-OPUC}. Without the three-term recurrence (\ref{3-term-recurrence}), one is forced to repeatedly compute the layer stripping iteration (\ref{layer-stripping}), which is computationally more expensive.   By contrast, implementation of the algorithm (\ref{sec-fast-algorithm}) in Matlab on an ordinary laptop takes 34 seconds to recover $c$ to near machine precision on a grid of $n=1000$ points.   

Theorem~\ref{thm-uniqueness} seems to be the first uniqueness result concerning inverse scattering of the Helmholtz equation for the class of piecewise constant wave speed using band limited data.  The near perfect performance of the inversion algorithm of \S\ref{sec-fast-algorithm} on numerical data with equal travel time layers confirms the paper's theoretical analysis.  What is surprising, however, is the efficacy of the same algorithm as an approximate inversion scheme for data which comes either from non-equal travel time layers or from continuously varying $c$.  This raises a basic theoretical question of how to account for the observed stability of the inversion, which is, moreover, reasonably robust in the presence of noise.  This is a topic for future investigation. The formulation of the inverse problem in the present paper is tailored to the practical scenario whereby data is collected at a single point for a finite range of frequencies.  A natural next step is to test the proposed inversion algorithm on experimental data.

\section*{Acknowledgments} 
The first two authors wish to thank the Dept.~of Mathematics \& Statistics at the University of Calgary for its generous hospitality during the summer of 2019. 



%
%
%
%
%
%

\end{document}